\documentclass[11pt]{amsart}
\newtheorem{theorem}{Theorem}[section]
\newtheorem{lemma}[theorem]{Lemma}
\newtheorem{proposition}[theorem]{Proposition}
\newtheorem{cor}[theorem]{Corollary}

\theoremstyle{definition}

\theoremstyle{remark}
\newtheorem{remark}[theorem]{Remark}

\numberwithin{equation}{section}

\newcommand{\abs}[1]{\lvert#1\rvert}

\newcommand{\supp}[1]{\operatorname{supp}{#1}}
\newcommand{\norm}[1]{\parallel\, #1\, \parallel}

\begin{document}

%\hfill{Preliminary version: Do not circulate please}

\title[Approximating L$^2$ invariants of amenable covering spaces]
{Approximating L$^2$ invariants of amenable covering spaces: A heat kernel approach.}
\author{Jozef Dodziuk}
\address{Ph.D.~Program in Mathematics, Graduate Center of CUNY, New York, NY 10036}
\email{jzdqc@cunyvm.cuny.edu}
\author{Varghese Mathai}
\address{Department of Mathematics, University of Adelaide, Adelaide 5005, Australia}
\date{JULY 1996}
\email{vmathai@maths.adelaide.edu.au}

\begin{abstract}
 In this paper, we prove that the $L^2$ Betti numbers of an amenable
covering space can be
approximated by the average Betti numbers of a regular exhaustion, under
some hypotheses. 
We also prove that some $L^2$ spectral invariants can be approximated 
by the corresponding average spectral invariants of a regular exhaustion.
The main tool which is used is a generalisation of the
"principle of not
feeling the boundary" (due to M. Kac), for heat kernels associated to
boundary value problems.

\end{abstract}

\maketitle

%------------------------------------------------------------
%------------------------------------------------------------

\section*{Introduction}

There has been important work done recently on approximating $L^2$ invariants on 
{\em residually finite} covering spaces due to L\"uck \cite{Lu} 
(see eg. \cite{Don1}
for earlier results). This paper is the first of
a couple of papers, where we study the approximation of $L^2$ invariants on {\em amenable} 
covering spaces, using heat kernel techniques. In a sequel paper \cite{DM}, we continue this 
study, but now using combinatorial techniques instead and where we manage to prove the 
main conjecture stated in this paper. 

Let $(M, g)$ be a compact Riemannian manifold and $\Gamma\rightarrow\widehat{M}\rightarrow M$
be a Galois covering space of $M$. Recall that {\em Cheeger's isoperimetric constant} is
defined as
$$
h(\widehat M, g) = \inf\Big\{\frac{{\mbox{vol}}_{n-1}(\partial D)}{{\mbox{vol}}_{n}(D)} : D \ \ {\text{ is open
in}}\ \  \widehat{M} \ \ {\text{and}} \ \ \partial D \ \ {\text{is smooth}}\Big\}.
$$

Now let $g'$ be another Riemannian metric on $M$. Since $g$ and $g'$ are quasi-isometric,
one sees that there is a positive constant $C$ such that 
$$
\frac{1}{C} h(\widehat M, g')\le h(\widehat M, g) \le C h(\widehat M, g').
$$

A Galois covering space $\Gamma\rightarrow\widehat{M}\rightarrow M$ is said to 
be {\em amenable} if $h(\widehat M, g) = 0$. By the previous discussion, this is
independent of the choice of the Riemannian metric $g$ on $M$. In fact, it turns out 
that $\widehat M$ is amenable if and only if $\Gamma$ is amenable.
There are many examples of amenable groups viz.

\begin{itemize}

\item[(1)] Abelian groups; 
\item[(2)] nilpotent groups and solvable groups;
\item[(3)] groups of subexponential growth;
\item[(4)] subgroups, quotient groups and extensions of amenable groups;
\item[(5)] direct limits of amenable groups.
\end{itemize}

A sequence of open, connected, relatively compact subsets $D_k, k=1,2,\ldots$ will
be called an {\em exhaustion} of $\widehat{M}$ if
\begin{itemize}
\item[(1)] $\overline{D}_k \subset D_{k+1}\hspace{5mm}\forall k \ge 1 $
\item[(2)] $\bigcup_{k\geq 1}D_k =\widehat{M}.$
\end{itemize}
 
If the group $\Gamma$ is amenable, we can choose an exhaustion such that
$ \lim_{k\rightarrow\infty}\frac{\mbox{vol}_{n-1}(\partial D_k)}
{\mbox{vol}_n (D_k)} = 0$ (cf. \cite{Ad}).  This is essentially a
combinatorial fact.  By choosing an exhaustion of the Cayley graph of
$\Gamma$ and a corresponding exhaustion of $\widehat M$ by unions of
translates of a fundamental domain with piecewise smooth boundary, and
smoothing the boundaries in a uniform way
 \cite[Proposition 3.2]{AdSu}, \cite{CG2}, we can assume that our exhaustion
has the following properties.

\begin{itemize}
\item[(1)]
For every fixed $\delta> 0$,  
\begin{equation}\label{bdry-layer} 
\lim_{k\rightarrow\infty}\frac{\mbox{vol}_{n}(\partial_\delta D_k)}
{\mbox{vol}_n (D_k)} = 0,
\end{equation}
where $\partial_\delta D_k = 
\{x\in\widehat{M} : d(x,\partial D_k) <\delta\}$.
\item[(2)]
The second fundamental forms of hypersurfaces $\partial D_k$ are uniformly
bounded by a constant independent of $k$.
\item[(3)]
The boundaries $\partial D_k$ are uniformly collared in the following sense.
There exists a constant $\rho$ such that the mapping
\begin{gather*}
\Phi_k : \partial D_k \times [0,\rho] \longrightarrow D_k\\
\Phi(x, t\xi) = \exp_x(t\xi)
\end{gather*}
is a diffeomorphism
onto its image.
 Here $x \in \partial D_k$, $ t \in
[0,\rho]$ and  $\xi$ is
the unit interior normal vector field on $\partial D_k$.
\end{itemize}

We shall call an exhaustion $\{ D_k \}^\infty_{k=1}$ satisfying these three
properties a {\em regular exhaustion} of $\widehat M$.

Let $b^j(D_k,\partial D_k)$ denote the $j^{\mbox{th}}$ Betti number of $D_k$
relative to the boundary and $b^j(D_k)$ denote the $j^{\mbox{th}}$ Betti number
of $D_k$.
Let $b_{(2)}^j(\widehat{M},\Gamma)$ denote the $j^{\mbox{th}}$ $L^2$ Betti
number
of $\widehat{M}$, then one of our main results is

\begin{theorem}
\label{thm:A}
Let $\Gamma\rightarrow\widehat{M}\rightarrow M$ be a noncompact
Galois cover of a compact manifold $M$ which is amenable.
Let $\{D_k\}_{k=1}^{\infty}$ be a regular exhaustion of $\widehat{M}$,
then for $j\geq 0$

\begin{itemize}
\item[a)]
$$
\limsup_{k\rightarrow\infty}\frac{\mbox{vol}(M)}
{\mbox{vol}(D_k)}b^j(D_k,\partial D_k)\leq b_{(2)}^j(\widehat{M},\Gamma).
$$
and
$$
\limsup_{k\rightarrow\infty}\frac{\mbox{vol}(M)}
{\mbox{vol}(D_k)}b^j(D_k)\leq b_{(2)}^j(\widehat{M},\Gamma).
$$
In particular, equality holds whenever $b_{(2)}^j(\widehat{M},\Gamma) = 0$.
\item[b)] More generally, one has for $N\geq 0$
$$
\limsup_{k\rightarrow\infty}\sum_{j=0}^N (-1)^{N-j} \frac{\mbox{vol}(M)}
{\mbox{vol}(D_k)}b^j(D_k,\partial D_k)\leq
\sum_{j=0}^N (-1)^{N-j} b_{(2)}^j(\widehat{M},\Gamma).
$$
and
$$
\limsup_{k\rightarrow\infty}\sum_{j=0}^N (-1)^{N-j} \frac{\mbox{vol}(M)}
{\mbox{vol}(D_k)}b^j(D_k)\leq
\sum_{j=0}^N (-1)^{N-j} b_{(2)}^j(\widehat{M},\Gamma).
$$
Moreover, there is equality when $N=n= \dim M$.
\end{itemize}
\end{theorem}

We remark that $b_{(2)}^j(\widehat{M},\Gamma) = 0$, trivially when $j=0,n$,
 but in many instances for other values of $j$ as well.

We conjecture that under the hypothesis of the theorem, one always has equality.

\vspace{.2in}

\noindent {\bf Conjecture.} {\em Let $\Gamma\rightarrow\widehat{M}\rightarrow M$ be a noncompact
Galois cover of a compact manifold $M$ which is amenable.
Let $\{D_k\}_{k=1}^{\infty}$ be a regular exhaustion of $\widehat{M}$,
then for $j\geq 0$
$$
\lim_{k\rightarrow\infty}\frac{\mbox{vol}(M)}{\mbox{vol}(D_k)}b^j
(D_k,\partial D_k)
=\lim_{k\rightarrow\infty}\frac{\mbox{vol}(M)}{\mbox{vol}(D_k)}b^j(D_k)
= b_{(2)}^j(\widehat{M},\Gamma ).
$$}

\vspace{.2in}

In other words, the conjecture says that on amenable covering spaces, 
the averaged absolute cohomology
and the averaged relative cohomology yield the same limit. 
The next theorem establishes this conjecture for manifolds of dimension
less than or equal to four, under a mild hypothesis.

\begin{theorem}
\label{thm:B}
Let $\Gamma\rightarrow\widehat{M}\rightarrow M$ be a noncompact Galois
cover of a compact manifold $M$ which is amenable and oriented.
\begin{itemize}
\item[a)] Suppose that $\mbox{dim} M = 2$. Then the conjecture is true. 

\item[b)] Suppose that $\mbox{dim} M = 3$ or $4$ and that ${\mbox{dim}} H^1(\widehat{M}) < \infty$.
Then the conjecture is true. 
\end{itemize}
\end{theorem}

\noindent{\bf Remarks.} There are many amenable Galois covering spaces of a 
surface of genus $\ge 1$. One sees this as follows. 
Since Galois covering spaces of a surface 
$\Sigma_g$ of of genus $g$ correspond to
normal quotient groups of the fundamental group $\pi_1(\Sigma_g) = \Gamma_g$,
the question of determining amenable Galois covering spaces of $\Sigma_g$ 
becomes a purely group theoretic question.
Some examples are as follows; one can consider the universal 
homology covering $\mathbb{Z}^{2g}\rightarrow\widehat{\Sigma_g}\rightarrow \Sigma_g$
and its quotients. A method of producing non-Abelian amenable covering
spaces is to first consider the natural homomorphism of the fundamental 
group $\Gamma_g$ onto the free group on $g$ generators $F_g$, and 
then one knows that there are many non-Abelian amenable groups on two
or more generators, i.e. quotients of $F_g$. For example $\mathbb{Z}/2 * \mathbb{Z}/2$, the integer 
Heisenberg groups etc. By the previous discussion, all these correspond
to non-Abelian amenable covering spaces of surfaces.  We remark that it
follows from Atiyah's $L^2$-index theorem \cite{A}
that $b^1_{(2)} \neq 0$ for
all these spaces.

We also establish the conjecture under the hypothesis that the heat kernel 
converges uniformly to the harmonic projection. More precisely, one has the following.

\begin{theorem}
\label{thm:C}
Let $\Gamma\rightarrow\widehat{M}\rightarrow M$ be a noncompact Galois
cover of a compact manifold $M$ which is amenable.
\begin{itemize}
\item[a)] Suppose that there is a function $f: [1,\infty] \to [0,\infty)$
which converges to zero as $t\to \infty$, such that for $j\ge 0$
$$
\frac{\mbox{vol}(M)}{\mbox{vol}(D_k)}\Big[\mbox{Tr}(e^{-t\Delta_j^{(k)}}) -
b^j(D_k,\partial D_k)\Big] \le f_j(t)
$$
for all $k\in\mathbb N$ and for all $t\ge 1$. Here $\Delta_j^{(k)}$ denotes 
the Laplacian on $D_k$ with the relative boundary conditions.
Then the conjecture is true.

\item[b)] Suppose that there is a function $f: [1,\infty] \to [0,\infty)$
which converges to zero as $t\to \infty$, such that for $j\ge 0$ 
$$
\frac{\mbox{vol}(M)}{\mbox{vol}(D_k)}\Big[\mbox{Tr}(e^{-t\Delta_j^{(k)}}) -
b^j(D_k)\Big] \le f_j(t)
$$
for all $k\in\mathbb N$ and for all $t\ge 1$. Here $\Delta_j^{(k)}$ denotes 
the Laplacian on $D_k$ with the absolute boundary conditions.
Then the conjecture is true.

\end{itemize}

\end{theorem}

In the first section, we recall some preliminary material on von Neumann
algebras and $L^2$ cohomology. In the second section, we prove the main
heat kernel estimates, which are quantitative versions of the principle 
of not feeling the boundary. In section 3, we 
give a heat kernel proof of Theorem \ref{thm:A}. We also give the proofs of
Theorems \ref{thm:B} and \ref{thm:C} in this section. Theorem \ref{thm:B}
uses Theorem \ref{thm:A} together with an index theorem. 
We also present non-examples to the analogue of Theorem \ref{thm:B} 
for non-amenable covering spaces.
A key result in the paper of Adachi and Sunada \cite{AdSu} is an approximation
property for the spectral density functions, with Dirichlet boundary conditions. They
also conjecture in section 1 of their paper that such an approximation property 
should also hold for the spectral density functions, with Neumann boundary conditions. 
Proposition 3.1 proves their conjecture, and much more, as complete results are 
proved for differential forms, and not just for functions.
In section 4, we derive
some useful corollaries of the Theorems which were proved in the earlier sections.
We also prove miscellaneous results on the spectrum
of the Laplacian on amenable covering spaces.

\section{The von Neumann trace and $L^2$ cohomology}

In this section,we briefly review results on the von Neumann trace and 
$L^2$ cohomology which are used in the paper. We shall assume 
that $\Gamma\rightarrow\widehat{M}\rightarrow M$ be a noncompact
Galois cover of a compact Riemannian manifold $(M,g)$.
Then the Riemannian metric
on $M$ lifts to a Riemannian metric on the Galois cover $\widehat M$, and
one can define the spaces of $L^2$ differential forms on $\widehat M$, \
$\Omega_{(2)}^\bullet (\widehat M)$. Let $d$ denote the deRham exterior 
derivative. One can define the space of {\em closed} $L^2$ differential
p-forms 
$$
Z^p(\widehat M) = \{\eta \in \Omega_{(2)}^p (\widehat M) : d\eta = 0\}.
$$
It is not hard to see that $Z^p(\widehat M) $ is a closed subspace of
$\Omega_{(2)}^p (\widehat M)$. One can also define 
the space of {\em exact} $L^2$ differential p-forms 
$$
B^p(\widehat M) = \{\eta \in \Omega_{(2)}^p (\widehat M) : \eta = d\phi\ \ {\text{for some}}
\ \ \phi\in \Omega_{(2)}^{p-1} (\widehat M)\}.
$$
In the example when $M$ is the circle, one can see that $B^p(\widehat M)$
is {\em not} a closed subspace of $\Omega_{(2)}^p (\widehat M)$. It is natural
then to take the closure, $\overline{B^p(\widehat M)}$, and to define the {\em reduced
$L^2$ cohomology} of $\widehat M$ as
$$
\bar{H}^p(\widehat M) =  Z^p(\widehat M)/\overline{B^p(\widehat M)}.
$$
It is clear that $\bar{H}^p(\widehat M)$ is $\Gamma$-invariant. 
Dodziuk \cite{D1} proved
that they depend only on the homotopy type of $M$.

Let $\delta$ denote the $L^2$ adjoint of $d$. Then the Laplacian
$\Delta_p$ acting on $\Omega_{(2)}^p (\widehat M)$ is
$$
\Delta_p = d\delta +\delta d.
$$
One can define the space of $L^2$ {\em harmonic forms} as
$$
{\mathcal H}^p(\widehat M) = \{\eta\in \Omega_{(2)}^p (\widehat M): \Delta_p\eta = 0\}.
$$
It is clear that ${\mathcal H}^p(\widehat M)$ is $\Gamma$-invariant and one can show that
${\mathcal H}^p(\widehat M)$ and $\bar{H}^p(\widehat M)$ are $\Gamma$-isomorphic Hilbert spaces.
One also has the Kodaira-Hodge decomposition (see \cite{GS})
$$
\Omega_{(2)}^p (\widehat M) = {\mathcal H}^p(\widehat M) \oplus 
\overline{d \Omega_{c}^{p-1} (\widehat M)}
\oplus \overline{\delta \Omega_{c}^{p+1} (\widehat M)}
$$
where $\Omega_{c}^\bullet (\widehat M)$ denotes the space of compactly supported
differential forms on $\widehat M$.

The commutant of ${\Gamma}$ action on 
$\Omega_{(2)}^p (\widehat M)$
is a von Neumann algebra with trace given
by
\[
\mbox{Tr}_{\Gamma}(L)=\int_{M}\mbox{tr}(L(x,x))d\mu(x),
\]
where $L$ is an operator in the commutant of $\Gamma$ 
with smooth
integral kernel $L(x,y)$, $x,y\in {\widehat M}$.
Therefore, there is a well defined dimension function 
$\mbox{dim}_{\Gamma}$, known as the {\em von Neumann dimension}
function,
which is defined on $\Gamma$-invariant closed subspaces
of $\Omega_{(2)}^p (\widehat M)$ and which takes 
values in $\mathbb R$. In particular, consider
the spectral projections $E_p(\lambda) = \chi_{[0,\lambda]}
(\Delta_p)$. By elliptic regularity theory, the integral
kernel of $E_p(\lambda)$ is smooth, and therefore one can define
the {\em von Neumann spectral function} for $\Delta_p$ as
$$
N_{p,\Gamma}(\lambda) := \mbox{dim}_{\Gamma}(\text{Im} E_p(\lambda))
= \mbox{Tr}_{\Gamma} (E_p(\lambda)) <\infty.
$$
In particular, the $L^2$ Betti numbers are defined as 
the von Neumann dimensions of the orthogonal projection onto the 
$L^2$ harmonic forms on $\widehat M$, $E_p(0)$
(cf. \cite{A}, \cite{D1}), i.e.
$$
b^p_{(2)}(\widehat M, \Gamma) = \dim_\Gamma({\mathcal H}^p(\widehat M))= 
\dim_\Gamma (\bar{H}^p(\widehat M)) 
= N_{p,\Gamma}(0).
$$
They have the following key properties.

\begin{itemize}
\item[a)] {\em Homotopy invariance property}, \cite{D1}.
For all $p\ge 0$,  $b^p_{(2)}(\widehat M, \Gamma)$ depends only on the homotopy type of 
the closed manifold $M$.
\item[b)] {\em Euler characteristic property}, \cite{A}.
$$\sum_{p=0}^n (-1)^p b^p_{(2)}(\widehat M, \Gamma) = \chi(M).$$
\item[c)] {\em Finite covers}.
If $\widehat M$ is a finite cover of $M$, then for all $p\ge 0$
$$b^p_{(2)}(\widehat M, \Gamma) = \frac{b^j(\widehat M)}{\# \Gamma}.$$
\item[d)] {\em Residually finite approximation}, \cite{Lu}.
Let $\ldots \Gamma_{m+1}\subset\Gamma_m\ldots \Gamma$ be a nested sequence
of normal subgroups of $\Gamma$ of finite index such that 
$\cap_{m\ge 1} \Gamma_m = \{e\}$. Then $\Gamma$ is said to be a {\em residually 
finite} group and one has for all $p\ge 0$
$$
b^p_{(2)}(\widehat M, \Gamma) = \lim_{m\to \infty}\frac{b^p(M_m)}{\#(\Gamma/\Gamma_m)}
$$
where $\Gamma/\Gamma_m \to M_m \to M$ is the corresponding sequence of finite
covers that "approximate" $\widehat M$.
\end{itemize}

\section{The principle of not feeling the boundary}

In this section, we shall assume 
that $\widehat{M}\rightarrow M$ is a noncompact Riemannian
covering of a compact manifold $M$ 
and $\{D_k\}_{k=1}^{\infty}$ is an exhaustion of $\widehat{M}$
satisfying conditions $(2)$ and $(3)$ in the definition of regular
exhaustion.  We emphasize that we do not require that the covering be
amenable, i.e. that (\ref{bdry-layer}) be satisfied. 
 We shall call the image $\Phi (\partial D_k \times [0,\rho])$ 
the collar of the boundary and
will denote it by $B(\partial D_k)$.
The key technical part of our paper is to prove 
a quantitative version of the following 

\vspace{.2in}

\noindent{\bf Principle of not feeling the boundary.} {\em
Let $k_j(t,x,y)$ denote the heat kernel on $L^2$ $j$-forms on $\widehat{M}$
and $p_j^k(t,x,y)$ denote the heat kernel on $D_k$ which is
associated to either the relative boundary conditions or the
absolute boundary conditions, then as $t\to 0$}
$$
k_j(t,x,y) \sim  p_j^k(t,x,y).
$$

\vspace{.2in}

This is a generalisation of the well known principle on functions due 
to M. Kac \cite{K}. The result below is weaker than what could be obtained
by the method of \cite[Section 7]{RS}.  The classical technique employed 
there 
is rather intricate and makes it difficult to identify the dependence of
constants in the estimates on geometry.  For applications that we have in
mind, it is essential to have estimates valid uniformly for all sets $D_k$
and the proof below does that.

\begin{theorem}{\bf (Main heat kernel estimate)}
\label{main-estimate}
Let $k_j(t,x,y)$ denote the heat kernel on $L^2$ $j$-forms on $\widehat{M}$
and $p_j^k(t,x,y)$ denote the heat kernel on $D_k$ which is
associated to either the relative boundary conditions or the
absolute boundary conditions. There are
positive constants $C_1,C_2$ depending only on the geometry
of $\widehat{M}$, the second fundamental form of $\partial D_k$, on the
choice of $\rho$ and on $T$, such that for $x,y \in D_k \setminus B(\partial D_k)$, $t\in (0,T]$,
\[
|k_j(t,x,y)-p_j^k(t,x,y)| \leq C_1e^{-C_2\frac{D^2(x)+D^2(y)}{t}},
\]
where $D_k(x)=d(x,\partial D_k)$.

\end{theorem}

To prove this theorem we study solutions of the heat equation in $D_k$
which vanish identically for $t=0$.  The difference of heat kernels
$
\omega (x,t)=k_j(t,x,y)-p_j^k(t,x,y)
$
for a fixed $y$  is such a solution as a consequence of
Minakshisundaram-Pleijel asymptotic expansion of the heat kernels
\cite{RS}.
The Weitzenb\"ock identity,$$
\Delta \omega = \nabla^{*}\nabla\omega + F_p\omega,$$
where $\nabla$ is the operator of covariant differentiation, $\nabla^{*}$
is its formal adjoint, and $F_p$ is an algebraic operator depending only
on the curvature of $M$, implies that 
$$
\left ( \Delta + \frac{\partial}{\partial t} \right )
e^{-\alpha t} | \omega |^2 =
e^{-\alpha t} (-2\langle \nabla \omega , \nabla \omega \rangle + 2 \langle
\omega , F_p \omega \rangle - \alpha |\omega |^2 ).
$$
Therefore, cf.\ \cite[Section 5]{D3}, for $\alpha = \sup_{x\in M} |F_p|$, 
$$
u(x,t) = e^{-\alpha t} |\omega (x,t)|^2
$$
satisfies
\begin{equation} \label{subsolution}
\left ( \Delta + \frac{\partial}{\partial t} \right ) u \leq 0
\end{equation}
i.e. $u(x,t)$ is a subsolution of the heat equation on functions.  In
particular, the weak maximum principle applies to $u$.  This is used in the
proof of the following lemma.

\begin{lemma}
\label{lemma-subsol}
Let $u(x,t)$ be a positive subsolution of the 
heat equation on functions in $D_k$
such that $u(x,0) \equiv 0$.  Suppose $$N = \sup_{(z,t)\in \partial D_{k}
\times [0,T]} u(z,t).$$  There exists constants $C_3,C_4 > 0$ 
independent of $k$
such that if $D(x) \geq \rho$ then $$
u(x,t) \leq C_3 N e^{-C_4 D^2(x)/t}$$
for $t\in [0,T]$.  The constants depend on $T$ and on the local geometry but
not on $k$.
\end{lemma}
\begin{proof}
Choose a smooth nonnegative function $\mu(r)$ on $[0,1]$  such that 
$\mu | [0,1/3]$ $\equiv 1$ and $\mu |[2/3,1] \equiv 0$.  For $x$ with 
$D(x) \geq \rho$ define a cut-off function $h(y)$ as follows. 
$$
h(y)=
\begin{cases}
0 & \text{if $d(x,y) \geq D(x)$},\\
1 & \text{if $d(x,y) \leq D(x) - (2/3)\rho$},\\
\mu(\frac{d(x,y) - D(x) + \rho}{\rho}) &
 \text{if $d(x,y) \in [D(x) - \rho, D(x)]$}.
\end{cases}
$$
Thus $h(y)$ vanishes near $\partial D_k$, and is identically equal to one
near $x$.  If $q(t,x,y)=k^0(t,x,y)$ is the heat kernel for functions on
$\widehat M$, using the fact that $u$ is a subsolution and that $q$ is a
solution of the heat equation, we see that
\begin{multline*} 
u(x,t)  = 
\int_0^t \frac{d}{d\tau}\int_{D_k}
h(y)u(y,\tau)q(t-\tau,y,x)\,d\mu(y)\,d\tau\\
  =  \int_0^t\int_{D_k}
\left ( h(y)\frac{du}{d\tau}(y,\tau) q(t-\tau,y,x)
+ h(y)u(y,\tau)\frac{dq}{d\tau}(t-\tau,y,x) \right )
\,d\mu(y)\,d\tau\\
 \leq \int_0^t\int_{D_k}
\left (-h(y)\Delta u(y,\tau) q(t-\tau,y,x) 
- h(y)u(y,\tau) \Delta_y q(t-\tau,y,x) \right )
\, d\mu(y)\,d\tau\\
  =  - 2\int_0^t\int_{D_k}
 u(y,\tau) \langle \nabla h(y),\nabla_yq(t-\tau,y,x)\rangle \,
d\mu(y)y\,d\tau\\
  \phantom{=} \qquad \qquad \qquad \qquad
- \int_0^t\int_{D_k}  u(y,\tau) \Delta h(y) q(t-\tau,y,x) 
\,d\mu(y)\,d\tau.
\end{multline*}
Note that according to the definition of $h(y)$ the integration on the
right hand side extends only
over the annullar region where $D(x)-(2/3)\rho \leq d(x,y) \leq D(x) 
-(1/3)\rho$.  In particular, $d(x,y) \geq D(x)/3 \geq \rho/3$ 
for such $y$.  Donnelly
proves ((see \cite{Don1}) that the heat kernel $q$ is almost euclidean, i.e.
that one has in particular
for all  $t\in (0,T]$
\begin{equation} \label{euclidean-don}
|\nabla_y^l\nabla_y^m q(t,x,y)| \leq 
C_5 t^{-(n+m+l)/2}e^{-\frac{C_6 d^2(x,y)}{t} } \qquad \mbox{for $k,l=0,1$}
\end{equation}
with constants $C_5,C_6$ depending only on the geometry of $M$ and $T$.
 We use this together with the maximum principle
applied to $u$ and the monotonicity of $e^{-C_6D^2(x)/s}$ 
as a function of $s$ 
to obtain
$$
u(x,t) \leq C_7 N e^{-C_8 D^2(x)/t} \mbox{vol}\,(B_{D(x)}(x)).
$$
Since the volume growth of $\widehat M$ is at most exponential, the lemma 
follows.
\end{proof}

To obtain Theorem \ref{main-estimate} from the lemma we need an estimate of the
form 
\begin{equation} 
\label{bound-bdry}
\sup_{(x,y,t)\in \partial D_k\times(D_k\setminus B(\partial D_k)\times[0,T]}
|k_j(t,x,y)-p_j^k(t,x,y)| \leq A,
\end{equation}
with $A$ independent of $k$.  We have such a bound for $k_j(t,x,y)$ in 
(\ref{euclidean-don}) so it suffices to estimate $p_j^k$.  We do this using the
finite propagation speed method as in \cite[Sections 1,2]{CGT}.  As observed by
Chernoff \cite{Chern}, solutions $\omega$ of the wave equation on $j$-forms
$$
\Delta_j \omega + \frac{\partial^2\omega}{\partial t^2} =0
$$ 
satisfying the initial condition 
$\frac{\partial\omega}{\partial t}(x,0) \equiv 0$ and
appropriate boundary conditions have finite
propagation speed i.e. if $\omega(x,0)$ is supported in $B_r(x)$, then
$\supp \omega(x,t) \subset B_{r+t}(x)$.  Chernoff considers only
complete manifolds with boundary but his proof carries over verbatim in
our context.  As explained in \cite{CGT}, under
appropriate bounded geometry assumptions, standard elliptic estimates (Sobolev
and G{\aa}rding inequalities) and finite propagation
speed give rise to pointwise estimates of
the heat kernel.  We review this briefly in our situation.  Consider
$\Delta=\Delta_j$
as an unbounded self-adjoint operator on $L^2$ $j$-forms on $D_k$ corresponding
to either relative or absolute boundary conditions.  By spectral theorem
applied to $\sqrt{\Delta}$ and
the Fourier inversion formula
\begin{align} \label{integral-formula}
\Delta^me^{-t\Delta}\Delta^l u & = \Delta^{m+l} e^{-t\Delta} u\\
{} & = \frac{1}{\pi}\int_0^\infty 
\frac{1}{\sqrt{2\pi t}} \frac{\partial^{2m+2l}}{\partial s^{2m+2l}}
\left ( e^{-\frac{s^2}{4t}}\right ) \cos (s\sqrt{\Delta}) u \,ds
\notag
\end{align}
for a fixed $t> 0$. Of course we require here that $u$ is in the domain
of $\sqrt{\Delta}$, i.e.\ it should satisfy appropriate boundary
conditions. Now
$\cos (s\sqrt \Delta)u$ satisfies the wave equation and the initial condition
above
so that, if $x_0 \in \partial D_k$, $y_0 \in
D_k\setminus B(\partial D_k)$ and $ \supp u \subset B_{\rho/3}(y_0)$,
$\supp \cos (t\sqrt\Delta )u \subset B_{\frac{\rho}{3}+t}(y_0)$.
The formula (\ref{integral-formula}) implies that
\begin{align*}
\norm{\Delta^m e^{-t\Delta} \Delta^l u}_{L^2(B_{\rho/3}(x_0)} & \leq
C_9 t^{-(\frac{1}{2} + 2m + 2l)} \norm{u} \int^\infty_{\rho/3}
e^{-\frac{s^2}{4t}}\,ds\\
{} & \leq
C_{10} t^{-(\frac{1}{2}+2m+2l)} e^{-\frac{\rho^2}{72 t}}
\leq C_{11} \norm{u}
\end{align*}
The kernel of $\Delta^m e^{-t\Delta} \Delta^l$ is $\Delta^m_x\Delta^l_y
p^k_j(t,x,y)$, so that if
$$
v(x,t)=\int_{D_k} \Delta^m p^k_j(t,x,y)u(y)\,d\mu(y),
$$
the inequality above gives an $L^2$ bound for $\Delta^m v$.  Summing over
$m=0,1\ldots [N/4]+1$ an applying standard elliptic estimates 
(which hold with uniform constants due to our bounded geometry assumptions)
we obtain an estimate for 
$\abs{v(x_0,t)}$ i.e.
$$
\abs{\int_{D_k} \Delta^l p_j^k(t,x_0,y)u(y)d\mu(y) } \leq C_{11} \norm{u}
$$
valid for smooth $u$ satisfying the boundary conditions and supported in
$B_{\rho/3}(y_0)$.  It follows that $$
\norm{\Delta^l p^k_j(t,x_0,y)}_{L^2({B_{\rho/3}(y_0)})} \leq C_{11}.
$$
Now, applying interior elliptic estimates,  we see again
that 
$$
\abs{ p_j^k(t,x_0,y_0)} \leq C_{12}.
$$ 
Since $x_0\in\partial D_k, y_0\in D_k \setminus B(\partial D_k)$ are arbitrary,
(\ref{bound-bdry}) is proven. 

\begin{proof}({\bf Theorem \ref{main-estimate}})
By (\ref{subsolution}), $u(x,t) =
e^{-\alpha t}\abs{p^k_j(t,x,y) -k_j(t,x,y)}^2$ with $y$ fixed satisfies 
the conditions of  Lemma \ref{lemma-subsol}.  Therefore
$$
\abs{p^k_j(t,x,y) -k_j(t,x,y)} \leq C_{1}e^{-C_{2}\frac{D^2(y)}{t}}
$$
by (\ref{bound-bdry}).
The theorem follows from the symmetry of the heat kernels in $x$ and $y$.
\end{proof}
In addition to the estimate above, we shall need a cheap bound of 
$\abs{p_j^k(t,y,y)}$
for $y$ near the boundary of $D_k$.  We actually state one valid for all
$y \in D_k$.
\begin{lemma}\label{diag-bound}
There exists a constant $c(t_0)$ which depends on the local geometry of $M$ but
is independent of $k$ such that for all $y \in D_k$ and $t \geq t_0$
$$
\abs{p_j^k(t,y,y)}\leq c(t_0). 
$$
\end{lemma}
\begin{proof}
Note that $x^{m+l}e^{-tx}\leq c(m,l)t^{-m-l}$ for all $x,t>0$.
Thus, by the spectral theorem,
$$
\norm{\Delta^m e^{-t\Delta}\Delta^l}\,\, \leq c(m,l)t^{-m-l}.
$$
This can be converted into a pointwise bound for the kernel exactly as above
and yields the lemma.
\end{proof}

\section{Proofs of the main theorems}

In this section, we shall also assume 
that $\Gamma\rightarrow\widehat{M}\rightarrow M$ is an noncompact
Galois cover of a compact manifold $M$ with an amenable covering group
$\Gamma$
and that $\{D_k\}_{k=1}^{\infty}$ is a regular exhaustion of $\widehat{M}$.

We begin with an approximation theorem for the heat kernels. A special case of
this proposition
settles a conjecture of Adachi and Sunada in \cite{AdSu}, section 1.

\begin{proposition} For every $t>0$, we have
\[
\mbox{Tr}_{\Gamma}(e^{-t\Delta_j}) =
\lim_{k\rightarrow \infty}\frac{\mbox{vol}(M)}{\mbox{vol}(D_k)}
\mbox{Tr}(e^{-t\Delta_j^{(k)}}),
\]
where $\Delta_j$ denotes the Laplacian acting on $L^2$ forms on
$\widehat{M}$ and $\Delta_j^{(k)}$ denotes the Laplacian on $D_k$ associated
to either the relative boundary conditions or the absolute 
boundary conditions.  The convergence is uniform in $t\in [t_0,t_1]$ 
for every $t_1>t_0>0$.
\end{proposition}

\begin{proof}
First observe that
\begin{eqnarray} \label{avg-ht}
\mbox{Tr}_{\Gamma}(e^{-t\Delta_j})&=&\int_F\mbox{tr}(k_j(t,x,x))d\mu (x)\\
&=& \lim_{k\rightarrow\infty}\frac{\mbox{vol}(M)}{\mbox{vol}(D_k)}
\int_{D_k}\mbox{tr}(k_j(t,x,x))d\mu (x)
\end{eqnarray}
where $F$ is a nice (say with piecewise smooth boundary) fundamental domain
for the action of $\Gamma$.
To see this we denote by $D'_k$ the union of those translates of
fundamental domains which are contained in $D_k$ and write 
\begin{align*}
{}& \mbox{Tr}_{\Gamma}(e^{-t\Delta_j})  = 
\frac{\mbox{vol}(M)}{\mbox{vol}(D'_k)}
\int_{D'_k}\mbox{tr}(k_j(t,x,x))d\mu (x)=\\
{}& \frac{\mbox{vol}(M)}{\mbox{vol}(D_k)}\int_{D_k} 
\mbox{tr}(k_j(t,x,x))d\mu (x) +
\frac{\mbox{vol}(M)}{\mbox{vol}(D_k)}\int_{D_k\setminus D'_k} 
\mbox{tr}(k_j(t,x,x))d\mu (x) \\
{}&\qquad \qquad
-\frac{\mbox{vol}(M)\mbox{vol}(D_k\setminus D'_k)}{
\mbox{vol}(D_k)\mbox{vol}(D'_k)}\int_{D'_k} \mbox{tr}(k_j(t,x,x))d\mu (x).
\end{align*}
Since $D_k \setminus D'_k\subset D_{\mbox{diam}(F)}$,
the last two summands tend to zero as $k\rightarrow\infty$ 
by (\ref{bdry-layer}) and (\ref{euclidean-don}).  This proves (\ref{avg-ht}).

We now estimate using Theorem \ref{main-estimate} and Lemma \ref{diag-bound},

\begin{eqnarray*}
& &\frac{\mbox{vol}(M)}{\mbox{vol}(D_k)}
\int_{D_k}\mbox{tr}(k_j(t,x,x)-p_j^k(t,x,x))d\mu (x)\\
&=& \frac{\mbox{vol}(M)}{\mbox{vol}(D_k)}
\int_{\partial_\delta D_k}\mbox{tr}(k_j(t,x,x)-p_j^k(t,x,x))d\mu (x)\\
&+& \frac{\mbox{vol}(M)}{\mbox{vol}(D_k)}
\int_{D_k\backslash \partial_\delta D_k}\mbox{tr}(k_j(t,x,x)
-p_j^k(t,x,x))d\mu (x)  \\
&\leq& c(t_0)\frac{\mbox{vol}(M)}{\mbox{vol}(D_k)}
\mbox{vol}(\partial_\delta D_k)
+ \frac{\mbox{vol}(M)}{\mbox{vol}(D_k)}C_{11}
e^{-\frac{C_{12}\delta^2}{t}}\mbox{vol}(D_k).
\end{eqnarray*}
Taking the limit as $k\rightarrow\infty$, we have
$$
\limsup_{k\rightarrow\infty}\big | 
\mbox{Tr}_{\Gamma}(e^{-t\Delta_j})
-\frac{\mbox{vol}(M)}{\mbox{vol}(D_k)}
\mbox{Tr}(e^{-t\Delta_j^{(k)}})
\big |
\leq 
\mbox{vol}(M) C_{11} e^{-\frac{C_{12}\delta^2}{t}}.
$$
for arbitrary $\delta > 0$, proving the proposition.
\end{proof}

\noindent{\bf Remarks.} A key result in the paper of Adachi and Sunada
 \cite{AdSu} 
is a special case of the approximation 
property in Proposition 3.1 above, with Dirichlet boundary conditions and for functions
(i.e. for the special case when $j=0$). They
also conjecture in section 1 of their paper that such an approximation property 
should also hold for the spectral density functions, with Neumann boundary conditions. 
Proposition 3.1 above proves their conjecture, and much more, as complete results are 
obtained for differential forms, and not just for functions.

\begin{proof}{(of Theorem \ref{thm:A})}
We will only prove (b), as (a) is obtained from (b) by adding the 2 inequalities
for $N=i$ and $N=i-1$.
Recall that $\Delta_j^{(k)}$ has a discrete point spectrum for all
$j\ge 0$, since the domain $D_k$ is compact.
Let $\lambda$ be a positive eigenvalue of $\Delta^{(k)}$, and $E_\lambda$
be the space
of eigenforms with eigenvalue $\lambda$. Then $E_\lambda$ decomposes into a
direct sum
$E_\lambda = E_\lambda^0 \oplus \cdots \oplus E_\lambda^n$, where
$E_\lambda^j$ denotes
the space of eigen $j$-forms of $\Delta_j^{(k)}$ with eigenvalue $\lambda$.
Since
$\lambda \ne 0$, the complex
$$
0\rightarrow E_\lambda^0\stackrel{d}\rightarrow \cdots \stackrel{d}\rightarrow E_\lambda^n \rightarrow 0
$$
is an {\em exact} sequence, so $D^N(\lambda) = \sum_{j=0}^N (-1)^{N-j} \dim
E_\lambda^j$
is non-negative for any $N$ such that $n\ge N\ge 0$. Also observe that
$$
\sum_{j=0}^N (-1)^{N-j}\mbox{Tr}(e^{-t\Delta_j^{(k)}}) =
\sum_{j=0}^N (-1)^{N-j}b^j(D_k,\partial D_k) + \sum_{\lambda\ne 0}
e^{-t\lambda}D^N(\lambda).
$$
Therefore we see that
$$
\sum_{j=0}^N (-1)^{N-j}\mbox{Tr}(e^{-t\Delta_j^{(k)}})\ge
\sum_{j=0}^N (-1)^{N-j}b^j(D_k,\partial D_k).
$$
Multiplying both sides of the previous inequality by $\frac{\mbox{vol}(M)}{\mbox{vol}(D_k)}$
and taking the limit superior as $k\to\infty$ and using Proposition 3.1, one sees that
$$
\sum_{j=0}^N (-1)^{N-j}\mbox{Tr}_{\Gamma}(e^{-t\Delta_j}) \ge
\limsup_{k\rightarrow\infty}\sum_{j=0}^N (-1)^{N-j}\frac{\mbox{vol}(M)}{\mbox{vol}(D_k)}b^j(D_k,\partial D_k).
$$
The proof of one of the inequalities part (b) is completed by taking the limit as 
$t\to \infty$.
The other inequality in part (b) is similarly proved.
\end{proof}

\begin{proof} {( of Theorem \ref{thm:B})}

We first prove the following ``index" theorems, (see \cite{J} for similar
results)
\begin{equation}
\lim_{k\rightarrow\infty}\frac{\mbox{vol}(M)}{\mbox{vol}(D_k)}
\chi(D_k, \partial D_k)=\chi(M)\label{eq:1}
\end{equation}
and
\begin{equation}
\lim_{k\rightarrow\infty}\frac{\mbox{vol}(M)}{\mbox{vol}(D_k)}
\chi(D_k)=\chi(M)\label{eq:2}
\end{equation}

These are proved as follows. 
Let $B$ denote either the absolute boundary conditions or the relative
boundary conditions, and $\chi(D_k)_B$ denote either $\chi(D_k)$ (for 
the absolute boundary conditions) and $\chi(D_k, \partial D_k)$ (for 
the relative boundary conditions). Then one has the McKean-Singer identity
(cf. \cite{G} section 4.2), for$t>0$,
\begin{equation}
\chi(D_k)_B=\sum_{j=0}^n (-1)^j \mbox{Tr}(e^{-t\Delta_j^{(k)}}).
\end{equation}
By Proposition 3.1, one sees that for $t>0$,
\begin{eqnarray*}
\lim_{k\rightarrow\infty}\frac{\mbox{vol}(M)}{\mbox{vol}(D_k)}\chi(D_k)_B
& = &\sum_{j=0}^n (-1)^j 
\lim_{k\rightarrow\infty}\frac{\mbox{vol}(M)}{\mbox{vol}(D_k)}
\mbox{Tr}(e^{-t\Delta_j^{(k)}})\\
& = & \sum_{j=0}^n (-1)^j \mbox{Tr}_{\Gamma}(e^{-t\Delta_j})\\
& = & \chi(M)
\end{eqnarray*}
where the last equality is the $L^2$ analogue of the McKean-Singer
identity (cf. \cite{Roe} chapter 13).
This establishes the equalities (3.3) and (3.4).
Note that when $n=\dim M$ is even, then $\chi(D_k)
= \chi(D_k, \partial D_k)$ and when $n$ is odd, then $\chi(D_k)
= - \chi(D_k, \partial D_k) = \frac{1}{2} \chi(\partial D_k)$.

Let $n=2$. Since $\Gamma$ is infinite
\[
b_{(2)}^0(\widehat{M},\Gamma)=0.
\]
By $L^2$ Poincar\'e duality (cf. \cite{A}, (6.4)) 
\[
b_{(2)}^2(\widehat{M},\Gamma)=0.
\]
So that by Theorem \ref{thm:A}, we see that Theorem \ref{thm:B} is true 
for $j\ne 1$. By (3.3) and (3.4), we see that
\begin{eqnarray*}
\lim_{k\rightarrow\infty}\frac{\mbox{vol}(M)}{\mbox{vol}(D_k)}
b^1(D_k,\partial D_k)&=&\chi (M) = \lim_{k\rightarrow\infty}\frac{\mbox{vol}(M)}{\mbox{vol}(D_k)}
b^1(D_k)\\
&=&b_{(2)}^1(\widehat{M},\Gamma)
\end{eqnarray*}
where we have used the Euler characteristic property of $L^2$ Betti numbers (cf. section 1).
This proves Theorem \ref{thm:B}, part a).

Let $n=3,4$. By a result of Cheeger and Gromov (cf. \cite{CG} Lemma 3.1), the natural
forgetful map
\[
H_{(2)}^1(\widehat{M})\hookrightarrow H^1(\widehat{M})
\]
is injective. Since $\Gamma$ is infinite, we know that $H_{(2)}^1(\widehat{M})$ is infinite
dimensional if $b_{(2)}^1(\widehat{M},\Gamma) \ne 0$. 
By hypothesis $\dim H^1(\widehat{M}) <\infty$, therefore
\[
b_{(2)}^1(\widehat{M},\Gamma)=0.
\]
Since $\Gamma$ is infinite
\[
b_{(2)}^0(\widehat{M},\Gamma)=0.
\]
Now let $n=3$.
By $L^2$ Poincar\'e duality (cf. \cite{A}, (6.4)) 
\[
b_{(2)}^3(\widehat{M},\Gamma)=0 = b_{(2)}^2(\widehat{M},\Gamma).
\]
So that by Theorem \ref{thm:A}, we see that Theorem \ref{thm:B} is true
when $n=3$.

Now let $n=4$. By $L^2$ Poincar\'e duality (cf. \cite{A}, (6.4)) 
\[
b_{(2)}^4(\widehat{M},\Gamma)=0 = b_{(2)}^3(\widehat{M},\Gamma).
\]
So that by Theorem \ref{thm:A}, we see that Theorem \ref{thm:B} is true 
for $j\ne 2$. By (3.3) and (3.4), we see that
\begin{eqnarray*}
\lim_{k\rightarrow\infty}\frac{\mbox{vol}(M)}{\mbox{vol}(D_k)}
b^2(D_k,\partial D_k)&=&\chi (M) = \lim_{k\rightarrow\infty}\frac{\mbox{vol}(M)}{\mbox{vol}(D_k)}
b^2(D_k)\\
&=&b_{(2)}^2(\widehat{M},\Gamma)
\end{eqnarray*}
where we have used the Euler characteristic property of $L^2$ Betti numbers (cf. section 1).
This proves Theorem \ref{thm:B}.

\end{proof}

We now give an example to show that the analogue of Theorem 0.2 is false 
(and hence so is the Conjecture) for certain non-amenable covering spaces.

\vspace{.2in}

\noindent{\bf Example.} 
Consider the $2n$-dimensional hyperbolic space $\mathbb{H}^{2n}$ and
the sequence of open balls $\{D_k\}_{k=1}^{\infty}$, where $D_k = B_k(x_0)$
denotes the ball of radius $k$ which is cetered at $x_0\in \mathbb{H}^{2n}$. 
Then, by an explicit calculation, 
$\{D_k\}_{k=1}^{\infty}$ is a regular exhaustion, i.e.
$$
\lim_{k\rightarrow\infty}\frac{\mbox{vol}_{n-1}(\partial D_k)}
{\mbox{vol}_n (D_k)} = h(\mathbb{H}^{2n}) = \frac{(2n-1)^2}{4}\ne 0.
$$
Let $\Gamma$ be a uniform lattice in $\mathbb{H}^{2n}$. It follows that
$\Gamma\to \mathbb{H}^{2n}\to M$ is a {\em non-amenable} Galois covering space, 
where $M= \mathbb{H}^{2n}/\Gamma$. It is a classical calculation 
that (see \cite{Don2}) that
$$
\mathcal{H}^j(\mathbb{H}^{2n}) = \left\{\begin{array}{lcl}
& & \{0\}\ \ \text{if} \ \ j\ne n;\\ \\
& & \text{infinite dimensional if} \ \ j=n. 
\end{array}\right.
$$
where $\mathcal{H}^j(\mathbb{H}^{2n})$ denotes the space of $L^2$
harmonic $j$-forms on $\mathbb{H}^{2n}$. It follows that
$$
b_{(2)}^n(\mathbb{H}^{2n}, \Gamma) \ne 0.
$$
However, since $D_k$ is topologically a ball in ${\mathbb R}^{2n}$, one
sees that for all $k\ge 1$,
$$
 H^n(D_k) = \{0\} = H^n(D_k, \partial D_k). 
$$
Therefore 
$$
\lim_{k\rightarrow\infty}\frac{\mbox{vol}(M)}{\mbox{vol}(D_k)}b^n(D_k)
= 0 < b_{(2)}^n(\mathbb{H}^{2n},\Gamma )
$$
and
$$
\lim_{k\rightarrow\infty}\frac{\mbox{vol}(M)}{\mbox{vol}(D_k)}b^n
(D_k,\partial D_k) = 0 < b_{(2)}^n(\mathbb{H}^{2n},\Gamma ).
$$\qed

\begin{proof} {( of Theorem \ref{thm:C})}
We first prove part a). The proof uses the hypothesis to justify the
interchange of the limits as $k\to\infty$ and as $t\to\infty$. We write 
$$
\Big| \frac{\mbox{vol}(M)}{\mbox{vol}(D_k)}
b^j(D_k,\partial D_k) - b_{(2)}^j(\widehat{M},\Gamma)\Big| 
$$
$$
\le f_j(t) + \Big|\mbox{Tr}_{\Gamma}(e^{-t\Delta_j})
-\frac{\mbox{vol}(M)}{\mbox{vol}(D_k)}
\mbox{Tr}(e^{-t\Delta_j^{(k)}})\Big| + 
\Big|\mbox{Tr}_{\Gamma}(e^{-t\Delta_j}) - b_{(2)}^j(\widehat{M},\Gamma)\Big|.
$$
Taking the limit as $k\to\infty$ and using Corollary \ref{avg-ht}, one has
$$
\Big| \lim_{k\rightarrow \infty}\frac{\mbox{vol}(M)}{\mbox{vol}(D_k)}
b^j(D_k,\partial D_k) - b_{(2)}^j(\widehat{M},\Gamma)\Big| 
\le f_j(t) +  
\Big|\mbox{Tr}_{\Gamma}(e^{-t\Delta_j}) - b_{(2)}^j(\widehat{M},\Gamma)\Big|.
$$
Since this is true for all $t\ge 1$, the Theorem follows by taking the limit
as $t\to \infty$ and using the hypothesis that $\lim_{t\to\infty} f_j(t) = 0$.
The proof of part b) is similar.

\end{proof}

\section{ Miscellaneous results on the spectrum of the Laplacian}

In this section, we derive some corollaries to the main theorems, as well
as several results on the spectrum of the Laplacian and zeta functions
associated to the Laplacian. We shall continue to assume 
that $\Gamma\rightarrow\widehat{M}\rightarrow M$ is an noncompact
Galois cover of a compact manifold $M$ with an amenable covering group
$\Gamma$ and that $\{D_k\}_{k=1}^{\infty}$ is a regular exhaustion of $\widehat{M}$.

Let $N_j^{(k)}(\lambda)$ denote the spectral function for $\Delta_j^{(k)}$, i.e.
$N_j^{(k)}(\lambda)$ is the number of eigenvalues of $\Delta_j^{(k)}$
which are less than or equal to $\lambda$. Let
$N_{j,\Gamma}(\lambda)$ denote the von Neumann's spectral function for
$\Delta_j$ (see section 1). Let
$\sigma(\Delta_j)$ denote the spectrum of $\Delta_j$ and 
$\sigma(\Delta_j^{(k)})$ denote the spectrum of $\Delta_j^{(k)}$. 
Then one has the following corollary of Proposition
\ref{avg-ht}.

\begin{cor} Let $\widehat M$ be an amenable Galois covering.
\begin{itemize}
\item[a)] At all points $\lambda$ which are points of continuity of 
$N_{j,\Gamma}(\lambda)$, one has
\[
\lim_{k\rightarrow\infty}\frac{\mbox{vol}(M)}{\mbox{vol}(D_k)}
N_j^{(k)}(\lambda)=N_{j,\Gamma}(\lambda).
\]

\item[b)] $\sigma(\Delta_j) \subset \overline{\bigcup_{k\ge 1} \sigma(\Delta_j^{(k)})}$
\end{itemize}
\end{cor}

\begin{proof} Part a) is immediately deduced from Proposition
\ref{avg-ht}, also using a Lemma due to Shubin (cf. \cite{Sh}). 
Now let $\lambda_1$ and $\lambda_2$ be points
of continuity of $N_{j,\Gamma}$, and $\lambda_1<\lambda_2$. Then by part a), one has
\[
\lim_{k\rightarrow\infty}\frac{\mbox{vol}(M)}{\mbox{vol}(D_k)}
\Big[N_j^{(k)}(\lambda_2) - N_j^{(k)}(\lambda_1)\Big]=N_{j,\Gamma}(\lambda_2)
- N_{j,\Gamma}(\lambda_1).
\]
From this, one easily deduces part b).
\end{proof}

\begin{remark} The limit $\lim_{k\rightarrow\infty}\frac{\mbox{vol}(M)}{\mbox{vol}(D_k)}
N_j^{(k)}(\lambda)$ is refered to in the Mathematical Physics literature as 
the {\em integrated density of states} for $\Delta_j$. The corollary above 
can be viewed as stating that the integrated density of states for $\Delta_j$
is independent of the choice of boundary conditions (relative or absolute).
\end{remark}
  
The closed subspaces 
$\overline{d \Omega_{c}^{p-1} (\widehat M)}$ and 
$\overline{\delta \Omega_{c}^{p+1} (\widehat M)}$
of $\Omega_{(2)}^p (\widehat M)$
are preseved by the Laplacian $\Delta_p$. Let 
$G_p(\lambda)$ and $F_p(\lambda)$
denote the respective spectral density functions, i.e.
$$
G_p(\lambda) = \mbox{Tr}_{\Gamma}(\chi_{[0,\lambda]}(\Delta_p^{<1>}))
$$ and 
$$
F_p(\lambda) = \mbox{Tr}_{\Gamma}(\chi_{[0,\lambda]}(\Delta_p^{<2>})) 
$$
where $\Delta_p^{<1>}$ denotes the retriction of $\Delta_p$ to the
closed subspace $\overline{d \Omega_{c}^{p-1} (\widehat M)}$ and 
$\Delta_p^{<2>}$ denotes the retriction of $\Delta_p$ on the
closed subspace $\overline{\delta \Omega_{c}^{p+1} (\widehat M)}$
and $\chi_{[0,\lambda]}$ denotes the characteristic function of the 
closed interval $[0,1]$. Then the following can be viewed as 
a generalization of one of the results in \cite{CG}. 

\begin{proposition}
Let $\widehat M$ be an amenable Galois covering 
and $\sigma(\Delta_p)$ denote the spectrum of the 
Laplacian acting on $L^2$ $p$-forms on $\widehat M$. Then $0\in \sigma(\Delta_1)$.
\end{proposition}

\begin{proof} We observe the following obvious characterization: 
$0 \in \sigma(\Delta_p)$ if and
only if $N_p(\lambda) >0$ for all $\lambda > 0$.
Since $\Gamma$ is amenable, one easily sees that $0\in \sigma(\Delta_0)$,
therefore $N_0(\lambda) >0$ for all $\lambda > 0$. 
By Lemma 3.1 in \cite{GS}, one
sees that $F_p(\lambda) = G_{p+1}(\lambda)$ for all $p\ge 0$.
So $N_1(\lambda) \ge G_{1}(\lambda) = F_0(\lambda) = N_0(\lambda) >0$,
i.e. $0\in \sigma(\Delta_1)$.  
\end{proof}

The following is a corollary of Theorem 0.3, noting that the large time
behaviour of the heat kernel corresponds to the small $\lambda$ behaviour 
of the corresponding spectral density function.

\begin{cor}Let $\widehat M$ be an amenable Galois covering. Suppose that there
is a positive constant $C$ such that for all $k\in\mathbb N$, one has
\[
\frac{\mbox{vol}(M)}{\mbox{vol}(D_k)}
\Big[N_j^{(k)}(\lambda) - N_j^{(k)}(0)\Big]\le C\lambda^{\beta_j}
\]
for all $\lambda \in (0,1)$. Then the $j^{th}$ Novikov-Shubin 
invariant (cf. \cite{GS}) \ $\alpha_j(\widehat M) \ge \beta_j$ and the conjecture
in the introduction is true.
\end{cor}

We now discuss zeta functions and convergence questions.

Let $\lambda >0$
and $\zeta(s, \Delta_p^{(k)} + \lambda)$ denote the zeta function of the 
operator $\Delta_p^{(k)} + \lambda$, upto normalization, that is,
$$
\zeta(s, \Delta_p^{(k)} + \lambda) = \frac{\mbox{vol}(M)}{\mbox{vol}(D_k)}
\frac{1}{\Gamma(s)}\int_0^\infty t^{s-1}\big(\mbox{Tr}
(e^{-t(\Delta_p^{(k)} + \lambda)}) - e^{-t\lambda} b^p(D_k,\partial D_k)\big) dt
$$
in the case of relative boundary conditions, and in the case of absolute
boundary conditions, one replaces $b^j(D_k,\partial D_k)$ by $b^j(D_k)$.
Then standard arguments \cite{RS} show that $\zeta(s, \Delta_p^{(k)} + \lambda)$
is a holomorphic function of $s$ in the half-plane $\Re(s)>n/2$, and it
has a meromorphic extension to $\mathbb C$ with no pole at $s=0$.

Also let $\lambda >0$
and $\zeta(s, \Delta_p + \lambda)$ denote the $L^2$ zeta function of
the operator $\Delta_p + \lambda$, that is,
$$
\zeta(s, \Delta_p + \lambda) = 
\frac{1}{\Gamma(s)}\int_0^\infty t^{s-1}\big(\mbox{Tr}_{\Gamma}
(e^{-t(\Delta_p + \lambda)}) - e^{-t\lambda} b^p_{(2)}(\widehat M)\big) dt.
$$
Then arguments as in \cite{M} show that $\zeta(s, \Delta_p + \lambda)$
is a holomorphic function of $s$ in the half-plane $\Re(s)>n/2$, and it
has a meromorphic extension to $\mathbb C$ with no pole at $s=0$.
Then the following proposition is relatively straightforward consequence 
of Proposition 3.1.

\begin{proposition}
For $\lambda>0$ fixed and as $k\to \infty$, $\zeta(s, \Delta_p^{(k)} + \lambda)$ 
converges to $\zeta(s, \Delta_p + \lambda)$, where the convergence is
uniform on compact subsets of $\Re(s)>n/2$.
\end{proposition}

If the following question can be answered in the affirmative, then one can
also deduce convergence of determinants as $k\to\infty$.

\noindent{\bf Question}. For $\lambda>0$ fixed and as $k\to \infty$, 
does $\zeta(s, \Delta_p^{(k)} + \lambda)$ 
converge to $\zeta(s, \Delta_p + \lambda)$, 
uniformly on compact subsets near $s=0$?

\end{document}